\documentclass{amsart}

\usepackage{mathrsfs}
\usepackage{amscd}
\usepackage{amsmath}
\usepackage{amssymb}
\usepackage{amsthm}
\usepackage{epsf}
\usepackage{latexsym}
\usepackage{verbatim}
\usepackage[all, cmtip]{xy}
\usepackage{tikz}
\usetikzlibrary{positioning}
\usetikzlibrary{matrix}
\usepackage{float}
\usepackage{hyperref}
\usepackage{comment}
\usepackage{bm}
\usepackage{ytableau}

\tikzstyle{bsq}=[rectangle, draw, thick, minimum width=.5cm, minimum height=.5cm]
\tikzstyle{bver}=[rectangle, draw, thick, minimum width=1cm, minimum height=2cm]
\tikzstyle{bhor}=[rectangle, draw, thick, minimum width=2cm, minimum height=1cm]

\usepackage[left=3.2cm, right=3.2cm]{geometry}

\newtheorem{theorem}{Theorem}[section]
\newtheorem{definition}[theorem]{Definition}

\newtheorem{conjecture}[theorem]{Conjecture}

\newtheorem{problem}[theorem]{Problem}
\newtheorem{varexample}[theorem]{Example}

\theoremstyle{definition}
\newtheorem{remark}[theorem]{Remark}

\def\bmu{\bm{\mu}}

\def\cM{\mathcal{M}}

\newcommand{\PP}{\mathbb{P}}
\newcommand{\M}{\overline{M}}

\newcommand{\cO}{\mathcal{O}}

\newcommand{\Pic}{\operatorname{Pic}}
\newcommand{\Sym}{\operatorname{Sym}}

\newenvironment{example}{\begin{varexample}
\begin{normalfont}}{\end{normalfont}
\end{varexample}}

\begin{document}
\title{Recent Developments in Brill-Noether Theory}
\author[D. Jensen]{David Jensen}
\address{David Jensen: Department of Mathematics,  University of Kentucky \hfill \newline\texttt{}
\indent 733 Patterson Office Tower,
Lexington, KY 40506--0027, USA}
\email{{\tt dave.jensen@uky.edu}}
\author[S. Payne]{Sam Payne}
\address{Sam Payne: Department of Mathematics,  University of Texas at Austin \hfill \newline\texttt{}
\indent 2515 Speedway, RLM 8.100,
Austin, TX 78712, USA}
\email{{\tt sampayne@utexas.edu}}
\date{October 18, 2021}
\bibliographystyle{alpha}

\maketitle

\begin{abstract}
We briefly survey recent results related to linear series on curves that are general in various moduli spaces, highlighting the interplay between algebraic geometry on a general curve and the combinatorics of its degenerations.  Breakthroughs include the  proof of the Maximal Rank Theorem, which determines the Hilbert function of the general linear series of given degree and rank on the general curve in $\mathcal{M}_g$, and complete analogs of the standard Brill-Noether theorems for curves that are general in Hurwitz spaces. Other advances include partial results in a similar direction for linear series in the Prym locus of a general unramified double cover of a general $k$-gonal curve and instances of the Strong Maximal Rank Conjecture.
\end{abstract}

\section{Linear Series on a Curve General in $\mathcal{M}_g$}
\label{Sec:BN}

In the early days of algebraic geometry, an algebraic curve was understood as coming equipped with an embedding in projective space.  In the late 19th and early 20th century, however, mathematicians began to consider abstract curves and Riemann surfaces, without reference to an ambient space.  Since that time, a central problem in algebraic geometry has been the following:

\begin{problem} \label{prob:BN}
Given a curve $C$, describe all maps from $C$ to projective space.
\end{problem}

When approaching this problem, it is natural to fix discrete invariants, such as the dimension of the projective space $\PP^r$ and the degree $d$ of the maps.  Passing to the language of line bundles and Picard varieties, we are then interested in the geometry of the \emph{Brill-Noether variety}
\[
W^r_d (C) \colon = \Big\{ L \in \Pic^d (C) \mbox{ } \vert \mbox{ } h^0 (C,L) \geq r+1 \Big\}. 
\]
A series of results in the 1980s describe these varieties when $C$ is general in the moduli space $\mathcal{M}_g$.  Throughout this survey, we assume that the ground field is algebraically closed field of characteristic zero. For extensions to positive characteristic, see Remark~\ref{rem:characteristic}.

\begin{theorem}
\label{Thm:BN}
Let $C$ be a curve that is general in $\mathcal{M}_g$.  Let $\rho (g,r,d) := g-(r+1)(g-d+r)$.
\begin{enumerate}
\item  The variety $W^r_d (C)$ is of pure dimension $\min\{\rho (g,r,d),g\}$ when this is non-negative, and empty otherwise \cite{GriffithsHarris80}{\em ;}
\item The locally closed variety $W^r_d (C) \smallsetminus W^{r+1}_d(C)$ is smooth \cite{Gieseker82}{\em ;}
\item  If $\rho (g,r,d) \geq 1$ then $W^r_d (C)$ is irreducible  \cite{FultonLazarsfeld81, Gieseker82}{\em ;}
\item  If $\rho(g,r,d) = 0$ then $\vert W^r_d (C) \vert = g! \cdot \prod_{j=0}^r \frac{j!}{(g-d+r+j)!}$  \cite{Kempf71, KleimanLaksov72, GriffithsHarris80}{\em ;}
\item If $\rho(g,r,d) = 0$ then the universal $W^r_d$ over the locus of curves in $\mathcal{M}_g$ with finitely many linear series of rank $r$ and degree $d$ is irreducible \cite{EisenbudHarris87}{\em ;}
%\item  If $\rho(g,r,d) < 0$ then every component of the locus in $\mathcal{M}_g$ of curves $C$ satisfying $W^r_d (C) \neq \emptyset$ has codimension at most $-\rho(g,r,d)$ \cite{Steffen98}. 
\end{enumerate}
\end{theorem}

This series of results about linear series on the generic curve in $\mathcal{M}_g$, and more recent generalizations and extensions to general curves in other moduli spaces, are broadly known as \emph{Brill-Noether Theory}, in honor of the celebrated work \cite{BN1874}. We will highlight the interplay between algebraic geometry and the combinatorics of degenerations in this subject. Note that the generic curve itself is a single object with no combinatorial structure, and it is possible to prove much of Theorem~\ref{Thm:BN} without degenerations, e.g., using vector bundle methods  \cite{Lazarsfeld86}.  Such an approach is used in \cite{ABFS16} to construct smooth Brill-Noether-Petri-general curves of every genus defined over $\mathbb{Q}$.    

In this brief survey, we emphasize the role of degenerations and their associated combinatorics in Brill-Noether theory.  The numerical invariants appearing in Theorem~\ref{Thm:BN} have rich combinatorial significance and intricate combinatorial structures come to the surface in various proofs by degeneration methods.  For instance, when $\rho(g,r,d) = 0$, the number of points in $W^r_d(C)$ is equal to the number of standard tableaux on the rectangular Young diagram $\lambda(g,r,d)$ with $r+1$ columns and $g-d+r$ rows.  A beautiful recent generalization says that the holomorphic Euler characteristic of $W^r_d (C)$ is equal to the number of set-valued tableaux on $\lambda(g,r,d)$ \cite{ACT21, CP21}.
%\begin{figure}[h!]
%\begin{ytableau}
%{} & {} \\
%{} & {} \\
%{} & {} \\
%\end{ytableau}
%\caption{The rectangular Young diagram $\lambda(6,1,4)$.}
%\end{figure} 

In general, the codimension of $W^r_d(C)$ in $\mathrm{Pic}^d(C)$ is the number of boxes in $\lambda(g,r,d)$. In some degeneration arguments, the limit of $W^r_d(C)$ breaks into pieces that are naturally in bijection with the standard fillings of $\lambda(g,r,d)$ with entries from $\{1, \ldots, g \}$. In others, the limiting pieces are related to solutions of an intersection problem on a Grassmannian, which is in turn related to the combinatorics of Young tableaux via Schubert calculus. The degeneration methods that have appeared in the Brill-Noether literature are themselves rich and varied. The limits of $C$ that have been fruitfully considered include $g$-nodal \cite{GriffithsHarris80, Gieseker82} and $g$-cuspidal rational curves \cite{EisenbudHarris83b}, flag curves with $g$ elliptic tails \cite{EisenbudHarris83c, EisenbudHarris87}, chains of elliptic curves \cite{Osserman14}, and tropical curves, most often chains of $g$ loops \cite{tropicalBN, tropicalGP}.  

Non-degenerative methods from algebraic geometry and the combinatorics of various  degenerations all continue to play key roles in current developments in Brill-Noether theory. In some cases, similarities in the combinatorial structures that appear through otherwise disjoint approaches are enough to yield new insights and breakthroughs. We hope readers will be enticed to dive more deeply into the details of this subject, to enjoy and appreciate the interplay between algebraic geometry and combinatorics in the study of special linear series on general curves.

\begin{remark} \label{rem:characteristic}
Here, we mention extensions of the surveyed results to positive characteristic. Theorems~\ref{Thm:BN} and~\ref{Thm:SMRC23} hold in arbitrary characteristic.  Theorem~\ref{Thm:PBN} holds in characteristic not 2.  Theorem~\ref{Thm:HBN} pre-supposes irreducibility of the Hurwitz space $\mathcal{H}_{k,g}$, which is known in characteristic greater than $k$ \cite{Fulton69}.  Nevertheless, the analog of parts (1)-(4) holds for some component of $\mathcal{H}_{k,g}$ in arbitrary characteristic, as does the analog of part (5) in any characteristic not dividing $k$.  Theorem~\ref{Thm:kPBN} pre-supposes irreducibility of the moduli of chains of covers; we are unaware of any extensions to positive characteristic.  Nevertheless, the analog of Theorem~\ref{Thm:kPBN} holds for some component of this moduli space in all characteristics not equal to $2$, $3$, or $5$ and not dividing $k$.  
\end{remark}

\noindent \textbf{Acknowledgments.} We are grateful to Gavril Farkas, Eric Larson, Hannah Larson, Yoav Len and Nathan Pflueger for helpful comments on an earlier draft of this survey.

The first author was supported by NSF DMS--2054135. The second author was supported by NSF DMS--2001502 and DMS--2053261.

\section{Hurwitz-Brill-Noether Theory}
\label{Sec:Hurwitz}

Theorem~\ref{Thm:BN} answers many natural questions about the geometry of $W^r_d(C)$ when $C$ is general in $\mathcal{M}_g$. But what about curves that are not general?  If a curve $C$ admits an atypical linear series, what does that imply about the existence and behavior of other linear series on $C$?  One compelling special case is the following analog of Problem~\ref{prob:BN} for curves with a special pencil.

\begin{problem}
\label{Prob:HBN}
If $\pi \colon C \to \PP^1$ is a general cover of genus $g$ and degree $k$, describe $W^r_d (C)$.
\end{problem}

\noindent Here, a general cover means that $C \to \PP^1$ is general in an appropriate Hurwitz space.  At first glance, all of the beautiful structures from classical Brill-Noether theory for curves general in $\mathcal{M}_g$ seem to break down.  One can, without much difficulty, construct examples where $W^r_d (C)$ is not equidimensional, and $W^r_d (C) \smallsetminus W^{r+1}_d(C)$ is singular.

\begin{example} \label{ex:trigonal}
Let $\pi \colon C \to \PP^1$ be a general cover of genus 6 and degree 3. Then $W^1_4 (C)$ is the disjoint union of a 1-dimensional component consisting of line bundles of the form $\pi^* \cO (1) \otimes \cO_C (p)$ for $p \in C$ and an isolated point $K_C \otimes \pi^* \cO (-2)$.
\end{example}

Nevertheless, special cases have been fully understood for many years, including all cases where $C \to \PP^1$ has degree $2$ or $3$ \cite{Clifford, Maroni}.  In the 1990s, researchers considered higher degree covers, found upper bounds for $\dim W^r_d(C)$, and proved the existence of components of certain specified dimensions \cite{Martens96, BallicoKeem96, CoppensMartens}.  More recently, Pflueger used tropical methods on a chain of loops to produce an improved upper bound for $\dim W^r_d(C)$, which he conjectured to be sharp  \cite{Pflueger17b}.  Pflueger also remarked that his results could be proved by similar combinatorial arguments using degeneration to a chain of elliptic curves instead of tropical degenerations, and both perspectives have been useful in subsequent work.

Pflueger's work inspired a new wave of interest in this area. His conjecture was proved by Ranganathan and the first author, using the algebraic geometry of curves in toric varieties in combination with tropical geometry and log deformation theory.

\begin{theorem}[\cite{JensenRanganathan}]
\label{Thm:JR}
Let $C \to \PP^1$ be a general cover of genus $g$ and degree $k$. Then
\[
\dim W^r_d (C) = \max_{\ell \in \{ 0, \ldots , r' \}} \rho (g,r-\ell,d) - \ell k ,
\]
where $r' = \min \{ r, g-d+r-1 \}$.  
\end{theorem}

\noindent Here, $W^r_d(C)$ may not have pure dimension, and $\dim W^r_d (C)$ means the maximum of the dimensions of its components.  This analog of Theorem~\ref{Thm:BN}(1) has been followed by more comprehensive results, giving a description of all components of $W^r_d(C)$, along with their dimensions and basic properties.  The starting point for this subsequent work is a stratification of $W^r_d (C)$, using the following discrete invariants that generalize the Maroni invariants of trigonal curves, introduced independently by H.~Larson \cite{Larson21} and by Cook-Powell and the first author \cite{CPJ19}. 

We say that a line bundle has \emph{splitting type} $\bmu = (\mu_1 , \ldots , \mu_k)$ if $\pi_* L \cong \oplus_{i=1}^k \cO (\mu_i)$.  The splitting type of a line bundle $L$ is a more refined invariant than its rank and degree; it encodes the rank of $L \otimes \pi^* \cO (m)$ for all integers $m$.  
\begin{definition}
The \emph{splitting type locus} is
\[
W^{\bmu} (C,\pi) \colon = \Big\{ L \in \Pic (C) \mbox{ } \vert \mbox{ } \pi_* L \cong \bigoplus_{i=1}^k \cO (\mu_i) \Big\} .
\]
\end{definition}

\begin{remark}
Note that $W^{\bmu}(C,\pi) \subset \Pic(C)$ is locally closed; it is the analog of $W^r_d(C) \smallsetminus W^{r+1}_d(C)$ in the setup for Theorem~\ref{Thm:BN}.  
\end{remark}

The analog of Theorem~\ref{Thm:BN} for splitting type loci involves numerical invariants that again have subtle and rich combinatorial significance, with $k$-core partitions playing the role of arbitrary partitions in standard Brill-Noether theory. While new in this context, $k$-core partitions have been studied extensively in number theory, representation theory, and combinatorics.  We recall the definition and a few basic properties; see \cite{LLMSSZ} for a fuller treatment and further references.

\begin{figure} [h!]
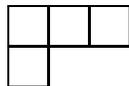

\begin{center}
\begin{ytableau}
{} & {} & {} \\
{} \\
\end{ytableau}
\end{center}
\caption{The Young diagram corresponding to the parition $4 = 3 + 1$.}
\label{Fig:4=3+1}
\end{figure}

We identify partitions with upper-left justified Young diagrams, so the partition $4 = 3 + 1$ is depicted by the Young diagram in Figure~\ref{Fig:4=3+1}.
The Young lattice $\mathcal{Y}$ is the set of partitions, partially ordered by inclusion of Young diagrams. 

Each box in a partition determines a \emph{hook}, consisting of that box, together with the boxes below it in the same column, and the boxes to the right of it in the same row.  The \emph{hook length} of a box is the number of boxes in the associated hook.  

\begin{definition}
A \emph{$k$-core partition} is a partition in which no box has hook length equal to $k$.
\end{definition}

\noindent For instance, the partition $4 = 3 + 1$ is $3$-core, but not $2$-core or $4$-core.

The $k$-core partitions $\mathcal{Y}_k \subset \mathcal{Y}$ form a sublattice with respect to the induced partial order. For a $k$-core partition $\lambda$, we write $r_k(\lambda)$ for the poset rank of $[\hat{0}, \lambda] \subset \mathcal{Y}_k$, i.e., the length of maximal chains in the interval between the empty partition $\hat{0}$ and $\lambda$.  Note that $r_k (\lambda)$ is typically strictly smaller than the number of boxes in $\lambda$, because deleting a corner box from a $k$-core partition typically does not yield another $k$-core partition.  

\begin{figure}[h!]
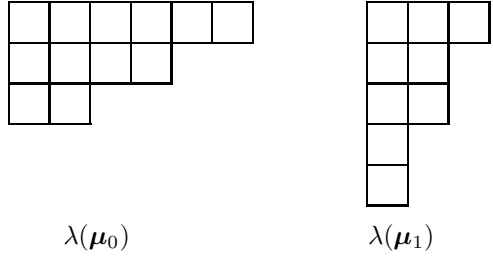

\begin{ytableau}
{} & {} & {} & {} & {} & {} \\
{} & {} & {} & {} \\
{} & {} \\
\end{ytableau}
\hspace{.5in}
\begin{ytableau}
{} & {} & {} \\
{} & {} \\
{} & {} \\
{} \\
{} \\
\end{ytableau} \\  \vspace{6 pt}
$\lambda(\bmu_0) \ \ \ \ \ \ \ \ \ \ \ \ \  \ \ \ \ \ \ \ \ \  \ \ \ \ \ \lambda(\bmu_1)$
\caption{The 3-core partitions associated to the splitting types $\bmu_0 = (-4,0,0)$ and $\bmu_1 = (-3,-2,1)$.}
\label{Fig:Trig}
\end{figure}

To each splitting type $\bmu = (\mu_1, \ldots, \mu_k)$, we associate a $k$-core partition $\lambda (\bmu)$ as follows:
\[
\lambda (\bmu) \colon = \Big\{ (x,y) \in \mathbb{N}^2 \mbox{ } \vert \mbox{ } \exists m \in \mathbb{Z} \text{ s.t. } x \leq \sum_{i=1}^k \max \{ 0, \mu_i + m + 1 \}, y \leq \sum_{i=1}^k \max \{ 0, -\mu_i - m - 1 \} \Big\} .
\]
For example, Figure~\ref{Fig:Trig} depicts two 3-core partitions of the form $\lambda (\bmu)$.

We can now state the main results on splitting type loci.

\begin{theorem} \label{Thm:HBN}
Let $\pi \colon C \to \PP^1$ be a general cover of genus $g$ and degree $k$.  Then
\begin{enumerate}
\item  The locally closed subvariety $W^{\bmu} (C,\pi) \subset \Pic^d (C)$ has pure dimension  $g - r_k (\lambda (\bmu))$ when this is nonnegative, and is empty otherwise \cite{Larson21, CPJ19}{\em ;}
\item If $g - r_k(\lambda(\bmu)) \geq 0$ then $W^{\bmu} (C,\pi)$ is smooth \cite{Larson21}{\em ;}
\item If $g - r_k(\lambda(\bmu)) \geq 1$ then $W^{\bmu} (C,\pi)$ is irreducible \cite{LLV}{\em ;}
\item If $g - r_k(\lambda(\bmu)) = 0$ then $\vert W^{\bmu} (C,\pi) \vert$ is the number of maximal chains in the interval $[\hat{0}, \lambda(\bmu)] \subset \mathcal{Y}_k$ \cite{LLV}{\em ;}
\item If $g - r_k(\lambda(\bmu)) = 0$ then the universal $W^{\bmu}$ over the locus of covers in the Hurwitz space with finitely many line bundles of splitting type $\bmu$ is irreducible \cite{LLV}.
\end{enumerate}
\end{theorem}

Theorem~\ref{Thm:HBN} is precisely analogous to Theorem~\ref{Thm:BN}, with the lattice of $k$-core partitions $\mathcal{Y}_k$ playing the role of the full lattice of partitions $\mathcal{Y}$.  Indeed, Theorem~\ref{Thm:BN} says that if $C$ is general in $\mathcal{M}_g$, then the codimension of $W^r_d (C)$ in $\Pic^d (C)$ is the rank of the interval $[\hat{0}, \lambda(g,r,d) ] \subset \mathcal{Y}$ and if $\dim W^r_d(C) = 0$ then $\vert W^r_d (C) \vert$ is equal to the number of maximal chains in this interval.  When $k$ is sufficiently large, Theorem~\ref{Thm:HBN} specializes to Theorem~\ref{Thm:BN}.  

\begin{example}
Returning to Example~\ref{ex:trigonal}, let $\pi \colon C \to \PP^1$ be a general cover of genus 6 and degree 3.  If $L \in W^1_4 (C)$ has splitting type $\bmu = (\mu_1 , \mu_2 , \mu_3)$, then
\[
\sum_{i=1}^3 \max \{ 0, \mu_i + 1 \} = h^0 (\pi_* L) = h^0 (L) = 2 \mbox{ and } 
\]
\[
\sum_{i=1}^3 \mu_i = \deg (\pi_* L) = \chi (L) -3 = -4 .
\]
The two splitting types that satisfy both these conditions and the condition that $6-r_3 (\lambda(\bmu)) \geq 0$ are $\bmu_0 = (-4,0,0)$ and $\bmu_1 = (-3,-2,1)$.  Each interval $[\hat{0}, \lambda(\bmu_i)] \subset \mathcal{Y}_3$ has a unique maximal chain, as shown in Figure~\ref{Fig:2Parts}. The lengths of these chains show that $r_3 (\lambda(\bmu_0)) = 6$ and $r_3 (\lambda(\bmu_1)) = 5$. By Theorem~\ref{Thm:HBN}, this implies that $\dim W^{\bmu_i} (C,\pi) = i$ and moreover, since the maximal chain in $[\hat{0}, \lambda(\bmu_0)] \subset \mathcal{Y}_3$ is unique, $W^{\bmu_0} (C,\pi)$ is a single point. A general trigonal curve is not hyperelliptic, so $W^2_4(C)$ is empty and thus, as discussed in Example~\ref{ex:trigonal}, $W^1_4 (C)$ has two irreducible components, one of dimension one and the other an isolated point of dimension 0.
\begin{figure}[H]

\begin{ytableau}
{1} & {2} & {3} & {4} & {5} & {6} \\
{3} & {4} & {5} & {6} \\
{5} & {6} \\
\end{ytableau}
\hspace{.5in}
\begin{ytableau}
{1} & {3} & {5} \\
{2} & {4} \\
{3} & {5} \\
{4} \\
{5} \\
\end{ytableau}

\caption{Maximal chains of 3-core partitions in the intervals $[\hat 0, \lambda (\bmu_0)]$ and $[\hat 0, \lambda (\bmu_1)]$.}
\label{Fig:2Parts}

\end{figure}
\end{example}

These results are proven by degeneration to a chain of elliptic curves in \cite{Larson21, LLV}.  To ensure that this degenerate curve admits a map of degree $k$ to a rational curve, one specifies that the difference between the two attaching points on each component has torsion order $k$ in the group law on the elliptic curve.  Building on Pflueger's earlier work, one can then classify the limit linear series on this degenerate curve that are limits of line bundles of a given splitting type.  The tableaux that arise from this construction are precisely those that come from maximal chains in the lattice of $k$-cores. The same combinatorial structures arise when approaching these problems via tropical degenerations, as in \cite{CPJ19, CPJ20}.

\begin{remark}
Another interesting result in late 20th century Brill-Noether Theory says that, if $\rho(g,r,d) < 0$ then every component of the locus in $\mathcal{M}_g$ of curves $C$ satisfying $W^r_d (C) \neq \emptyset$ has codimension at most $-\rho(g,r,d)$ \cite{Steffen98}. To the best of our understanding, the analogs of this statement in Hurwitz-Brill-Noether Theory, and in the Prym-Brill-Noether theory discussed below, are open problems.  
\end{remark}

\section{Prym-Brill-Noether Theory}
\label{Sec:Prym}

Let $C$ be a smooth curve of genus $g$, and let $\pi \colon \widetilde{C} \to C$ be an unramified double cover.  In \cite{Welters85}, Welters defines the \emph{Prym-Brill-Noether locus} to be
\[
V^r (C,\pi) \colon = \Big\{  L \in \Pic^{2g-2} (C) \mbox{ } \vert \mbox{ } \textrm{Nm}_{\pi} (L) = \omega_C , h^0 (C,L) \geq r+1, \text{ and } h^0 (C,L) \equiv r+1 \pmod{2} \Big\} .
\]
The following is an analog of parts (1)-(4) of Theorem~\ref{Thm:BN} in this setting.

\begin{theorem}
\label{Thm:PBN}
Let $\pi \colon \widetilde{C} \to C$ be a general unramified double cover, where the base has genus $g$.
\begin{enumerate}
\item  The variety $V^r (C,\pi)$ is of pure dimension $g-1 - {{r+1}\choose{2}}$ when this is non-negative, and empty otherwise \cite{Welters85, Bertram87}{\em ;}
\item  The singular locus of $V^r (C,\pi)$ is $V^{r+2} (C,\pi)$ \cite{Welters85} {\em ;}
\item  If $g-1 - {{r+1}\choose{2}} \geq 1$ then $V^r (C,\pi)$ is irreducible \cite{Debarre00} {\em ;}
\item  If $g-1 - {{r+1}\choose{2}} = 0$, then $\vert V^r (C,\pi) \vert = 2^{{r}\choose{2}} \cdot (g-1)! \cdot \prod_{i=1}^r \frac{(i-1)!}{(2i-1)!}$ \cite{dCP95}.
\end{enumerate}
\end{theorem}

\noindent It is possible to prove (3) and (4) using vector bundle techniques rather than degenerations \cite{dCP95}. As in Theorem~\ref{Thm:BN}, the numerical invariants that appear have rich combinatorial significance, and intricate structures involving partitions and tableaux emerge in proofs via degeneration.  

\begin{remark}
The Prym-Brill-Noether analog of Theorem~\ref{Thm:BN} (5) is an open problem.
\end{remark}

Let $\Delta_r = r + (r-1) + \cdots + 1$ denote the isosceles triangular partition with side length $r$.
\begin{figure}[h!]
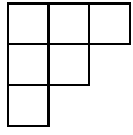


\begin{ytableau}
{} & {} & {} \\
{} & {} \\
{} \\
\end{ytableau}

\caption{The partition $\Delta_3$.}
\label{Fig:Prym}

\end{figure}
This partition plays a role in Prym-Brill-Noether theory analogous to that of the rectangular partition $\lambda(g,r,d)$ in the Brill-Noether theory of a general curve in $\cM_g$.  Indeed, Theorem~\ref{Thm:PBN} says that, if $\pi \colon \widetilde{C} \to C$ is a general unramified double cover, then the codimension of $V^r (C,\pi)$ in the Prym variety is the rank of the interval $[\hat{0}, \Delta_r] \subset \mathcal{Y}$, and if $\dim V^r (C,\pi) = 0$, then $\vert V^r (C,\pi) \vert$ is equal to the number of maximal chains in this interval.

Given Theorem~\ref{Thm:PBN} for unramified double covers of a general curve, it is natural to ask next about the Prym-Brill-Noether theory of double covers of a curve $C$ that is not necessarily general in $\cM_g$. As in Section~\ref{Sec:Hurwitz}, a compelling special case is when $C$ is general in the $k$-gonal locus, i.e., the irreducible locus of curves in $\cM_g$ on which the minimal degree of a pencil is $k$, for $k \geq 3$.  When $k \geq 3$, there is an irreducible moduli space of chains of covers $\widetilde{C} \to C \to \PP^1$, where the first map is an unramified double cover and the second has degree $k$ \cite{BF86}, so the generic such cover is well-defined.  %The closure of the image of this space in the moduli space of double covers $\mathcal{R}_g$ is known as the \emph{$k$-gonal locus} in $\mathcal{R}_g$.  
%For covers in the $k$-gonal locus, 
In this setting, recent work provides new bounds on the dimension of $V^r (C,\pi)$, obtained for even $k$ in \cite{LenUlirsch21} and for all $k$ in \cite{CLRW20}.

\begin{theorem} \cite{LenUlirsch21, CLRW20}
\label{Thm:kPBN}
Let $\pi \colon \widetilde{C} \to C$ be an unramified double cover of a general $k$-gonal curve of genus $g$, with $k \geq 3$.  Let $\ell = \left \lceil \frac{k}{2} \right \rceil$, and let 
\[
n(r,k) \colon = \left\{ \begin{array}{ll}
{{\ell+1}\choose{2}} + \ell (r-\ell) & \textrm{if $\ell \leq r-1$} \\
{{r+1}\choose{2}} & \textrm{if $\ell \geq r$.}
\end{array} \right.
\]
Then
\[
\dim V^r (C, \pi) \leq g-1-n(r,k) .
\]
\end{theorem}

\noindent This upper bound is not  sharp in general.  For example, consider the case where $g=6$ and $k=4$.  A general curve of genus 6 is tetragonal, so by Theorem~\ref{Thm:PBN} we see that $V^3 (C,\pi)$ is empty.  Since $\ell = 2 = r-1$, however, Theorem~\ref{Thm:kPBN} yields only
\[
\dim V^3 (C, \pi) \leq 6-1-{{3}\choose{2}} - 2 = 0 .
\]
However, the bound is conjectured to be sharp when $g$ is sufficiently large.

\begin{conjecture} \cite{CLRW20} %\cite[Conjecture~3.9]{CLRW20}
Let $\pi \colon \widetilde{C} \to C$ be an unramified double cover of a general $k$-gonal curve, with $k \geq 3$.  If $g \gg n(r,k)$, then 
\[
\dim V^r (C, \pi) = g-1-n(r,k).
\]
\end{conjecture}

\begin{figure}[H]
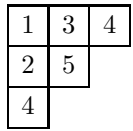


\begin{ytableau}
{1} & {3} & {4} \\
{2} & {5} \\
{4} \\
\end{ytableau}

\caption{A maximal chain in $\mathcal{Y}_4$.}
\label{Fig:kPrym}

\end{figure}

From a combinatorial perspective, Theorem~\ref{Thm:kPBN} arises from studying the partition $\Delta_r$ in the lattice of $k$-core partitions $\mathcal{Y}_k$.  Note that every box in $\Delta_r$ has odd hook length, so if $k$ is even, then $\Delta_r$ is a $k$-core partition.  When $k$ is even, $n(r,k)$ is the rank of $\Delta_r$ in $\mathcal{Y}_k$.   For example, Figure~\ref{Fig:kPrym} depicts a maximal chain in the interval $[ \hat{0}, \Delta_3] \subset \mathcal{Y}_4$.  This chain has length 5, corresponding to the fact that $n(3,4) = r_4 (\Delta_3) = 5$.

A similar result holds when $k$ is odd, but because $\Delta_r$ is not a $k$-core in this case, it requires the language of displacement tableaux from \cite{Pflueger17b}.  A tableau $t$ is said to be a \emph{$k$-uniform displacement tableau} if, whenever $t(x,y) = t(x',y')$, we have $x-y \equiv x'-y' \pmod{k}$.  The value $n(r,k)$ is equal to the minimum number of symbols in a $k$-uniform displacement tableau on $\Delta_r$.

Theorem~\ref{Thm:kPBN} is proven via tropical methods.  In \cite{LenUlirsch21}, Len and Ulirsch show that the skeleton of a Prym variety is the Prym variety of the skeleton.  They then construct an explicit double cover $\varphi \colon \widetilde{\Gamma} \to \Gamma$ of tropical curves, which they call the \emph{$k$-gonal uniform folded chain of loops}, and study the tropical Prym-Brill-Noether variety $V^r (\Gamma, \varphi )$.  These ideas are developed further in \cite{CLRW20}, which establishes not only the dimension bound above, but also various topological properties.% of the tropical Prym-Brill-Noether variety.

\begin{theorem} \cite{CLRW20}
\label{Thm:PrymConnect}
Let $\varphi \colon \widetilde{\Gamma} \to \Gamma$ be a $k$-gonal uniform folded chain of loops.  Then the tropical Prym-Brill-Noether variety $V^r (\Gamma, \varphi)$ has pure dimension $g-1-n(r,k)$ and, when this dimension is positive, it is connected in codimension one.
\end{theorem}

\noindent Tropicalizations of irreducible varieties are of pure dimension and connected through codimension 1 \cite{CartwrightPayne12}, so Theorem~\ref{Thm:PrymConnect} points suggestively toward analogs of parts (1) and (3) of Theorem~\ref{Thm:BN}.% in this context.

Many open questions remain about the Prym-Brill-Noether variety $V^r (C,\pi)$, where $\pi \colon \widetilde{C} \to C$ is an unramified double cover of a general $k$-gonal curve.

\begin{problem}
Let $\pi \colon \widetilde{C} \to C$ be an unramified double cover of a general $k$-gonal curve, with $k \geq 3$.
\begin{enumerate}
\item  What is $\dim V^r (C,\pi)$?
\item  Is the locally closed variety $V^r (C,\pi) \smallsetminus V^{r+2} (C,\pi)$ smooth?
\item  If $\dim V^r (C,\pi) \geq 1$, is $V^r (C,\pi)$ irreducible?
\item  When $\dim V^r (C,\pi) = 0$, what is $\vert V^r (C,\pi) \vert$?
\item  When $\dim V^r (C,\pi) = 0$, what is the action of monodromy?
%\item When the expected dimension of $V^r(C, \pi)$ is negative, what are the dimensions of the components of the locus of pairs $(C', \pi')$ such that $V^r(C', \pi')$ is not empty?
\end{enumerate}
\end{problem}

One might reasonably approach such problems starting from \cite{CLRW20} via tropical lifting, as in \cite{LiftingDivisors}, or via log deformation theory, as in \cite{JensenRanganathan}.  Alternatively, one could look for analogous statements about limit linear series on a double cover of a $k$-gonal chain of elliptic curves by a folded chain of elliptic curves, and identify which of these are limits of line bundles in the Prym-Brill-Noether locus of a degenerating family of covers, as was done in \cite{LLV} for the Hurwitz-Brill-Noether theory.

\section{Multiplication Maps and Maximal Rank Results}
\label{Sec:MRC}

A \emph{linear series} of rank $r$ on a curve $C$ is a pair $(L,V)$, where $L$ is a line bundle on $C$ and $V \subseteq H^0(C,L)$ is a vector subspace of dimension $r+1$.  By Theorem~\ref{Thm:BN}, when $\rho(g,r,d) \geq 0$, the universal space $G^r_d$ of linear series of degree $d$ and rank $r$ has a unique irreducible component that dominates $\mathcal{M}_g$. In other words, there is a well-defined generic linear series on the generic curve of genus $g$. Natural directions for further study therefore include the syzygies of this generic linear series, and the geometry of the locus of linear series with special syzygies. 

The most basic syzygies to consider in this context are the kernels of the multiplication maps from symmetric powers of a linear series to global sections of tensor powers of the underlying line bundle.  Recall that a map between finite dimensional vector spaces has \emph{maximal rank} if it is either injective or surjective.  Harris popularized the \emph{Maximal Rank Conjecture}, i.e., the prediction that the multiplication maps for the general linear series on the generic curve should have maximal rank \cite{Harris82}. This is now a theorem of E.~Larson.

\begin{theorem} [\cite{Larson17}] \label{Thm:MRC}
Suppose $C$ is general and $V \subset H^0(C,L)$ is a general linear series of given degree and rank.  Then the multiplication maps
\[
\mu_m: \Sym^m V \rightarrow H^0(C,L^{\otimes m})
\]
have maximal rank for all $m$.
\end{theorem}

The proof of this theorem relies on an elaborate induction involving interpolation for vector bundles \cite{AtanasovLarsonYang19} and an analogous statement for hyperplane sections of the general curve of genus $g$ and degree $d$ in $\PP^r$ \cite{Larson20}. It also uses novel degeneration techniques that surpass what was previously achieved using embedded methods \`a la Hirschowitz-Horace \cite{Hirschowitz85}, in part by exploiting properties of the Kontsevich moduli space of stable maps. The essential ingredients span more than a half-dozen articles, some of which remain unpublished, as does the final paper completing the proof.  For a short and helpful overview of this work, see \cite{Larson18}.

When $\rho(g,r,d) = 0$, Theorems~\ref{Thm:BN} and \ref{Thm:MRC} together imply that every linear series of degree $d$ and rank $r$ on a general curve has multiplication maps of maximal rank for all $m$.  When $\rho(g,r,d)$ is positive, other natural questions arise about the geometry of the locus in $G^r_d(C)$, for a general curve $C$, where a given multiplication map fails to have maximal rank. In this direction, Aprodu and Farkas proposed a \emph{Strong Maximal Rank Conjecture} that predicts the codimension of the locus in $G^r_d(C)$ where $\mu_m$ fails to have maximal rank, when $C$ is general, for suitable parameters $g,r,d,m$.

\begin{conjecture} [\cite{AproduFarkas11}]
Fix integers $g,r,d \geq 1$ such that $0 \leq \rho(g,r,d) < r-2$, and fix $m \geq 2$.  Let $C$ be a general curve in $\mathcal{M}_g$.  Then the determinantal variety 
\[
\{ (L,V) \in G^r_d(C) \mid \mu_m \colon \Sym^m V \to H^0(C, L^{\otimes m}) \text{  fails to have maximal rank } \}
\]
is of the expected dimension
\[
\rho(g,r,d) - 1 - \left|{r+m \choose m} - (2d+1-g) \right|,
\]
where, by convention, the locus is empty when the expected dimension is negative.
\end{conjecture}

One motivation for pursuing this conjecture is that the images of the loci of pairs $(C,L)$ for which $\mu_m$ fails to have maximal rank are candidates for interesting effective cycles in the coarse moduli space $\M_g$. When these syzygytic loci are of the expected dimension, their classes can sometimes be computed via intersection theory and computations with test curves. For instance, the counterexamples to the Slope Conjecture in \cite{FarkasPopa05, Farkas09b} are instances of such syzygytic loci that happen to be effective divisors outside the cone spanned by previously known effective divisors.

The cases of the Strong Maximal Rank Conjecture where $\rho(g,r,d) = 0$ follow from Theorem~\ref{Thm:MRC}. Only a small selection of other cases are known; these include the non-special cases where $r = d-g$, \cite[Proposition~5.7]{AproduFarkas11}, and the cases where $m = 2$ and $d  \leq g+1$ \cite{TiB03}.   Note that, for $m = 2$, and for each fixed $r$, there are only finitely many other cases to consider. The cases where $m=2$ and $r=4$ are discussed in \cite[Remark~5.6]{AproduFarkas11}. 

Farkas and the authors recently proved a few new cases of the Strong Maximal Rank Conjecture for $m = 2$, starting with $(g,r,d) = (22,6,25)$ and $(23,6,26)$. In these cases, every linear series of the given parameters on a general curve $C$ of genus $g$ is complete, so $W^r_d(C) = G^r_d(C)$, and the expected dimension of the locus where the multiplication map fails to have maximal rank is negative.

\begin{theorem} [\cite{M23}]
\label{Thm:SMRC23}
Let $g=22$ or $23$. For a general curve $C$ in $\cM_g$, the multiplication map
\[
\mu_2 \colon \Sym^2 H^0(C, L) \to H^0(C,L^{\otimes 2})
\]
is injective for \emph{all} line bundles $L\in W^6_{g+3}(C)$.
\end{theorem}

Once again, the syzygytic loci obtained as images of pairs $(C,L)$ where the multiplication map fails to be injective turn out to be effective divisors of smaller slope than any previously known effective divisors on $\M_{22}$ and $\M_{23}$, respectively. Moreover, the slopes are small enough to conclude that these moduli spaces are of general type.

\begin{theorem} [\cite{M23}]
\label{Thm:GeneralType}
The moduli space of curves $\M_g$ is of general type for $g \geq 22$.
\end{theorem}

\noindent This improves the bound of $g \geq 24$ proved decades ago by Harris, Mumford, and Eisenbud \cite{HarrisMumford82, Harris84, EisenbudHarris87b}; they used Theorem~\ref{Thm:BN} to construct the requisite divisors of small slope in those cases. The case $g = 22$ also disproves the conjecture of Eisenbud and Harris that $\M_g$ should be uniruled for $g < 23$ \cite{EisenbudHarris89b}.  The Kodaira dimensions and uniruledness of $\M_g$ for $17 \leq g \leq 21$ remain open problems. 

\medskip

Farkas and the authors also applied the techniques from \cite{M23}, with some small but important improvements, to prove the Strong Maximal Rank Conjecture for $m =2$ and $(g,r,d) = (13,5,16)$, and thereby construct divisors of small slope on $\M_{13}$.

\begin{theorem}[\cite{R13}]
For a general curve $C$ in $\cM_{13}$, the multiplication map
\[
\mu_2 \colon \Sym^2 H^0(C, L) \to H^0(C,L^{\otimes 2})
\]
is surjective for all line bundles $L \in W^5_{16}(C)$.
\end{theorem}

The closure of the locus of Brill-Noether general curves of genus 13 with at least one $L \in W^5_{16}$ such that $\mu_2$ is not surjective supports an effective divisor in $\M_g$.  This is the first known effective divisor on a moduli space $\M_g$ for $g \geq 10$ with slope less than $6 + 10/g$; it had previously been believed that such divisors may not exist \cite{ChenFarkasMorrison13}.  The class of this divisor is used to compute the class of yet another divisor of even smaller slope, the \emph{non-abelian Brill-Noether divisor}, which is used to show that the  Prym moduli space $\overline{\mathcal{R}}_{13}$ is of general type; 13 is now the smallest $g$ for which $\overline{\mathcal{R}}_g$ is known to be of general type.   

\medskip

This survey is far from exhaustive, and highlights only a portion of the recent results in this active area.  It inevitably reflects the biases and limitations of the authors. We hope the references provided will be a useful starting point for those who wish to learn more.

\bibliography{BN-RecentDevs}

\newcommand{\etalchar}[1]{$^{#1}$}
\begin{thebibliography}{CLRW20}

\bibitem[ABFS16]{ABFS16}
E.~Arbarello, A.~Bruno, G.~Farkas, and G.~Sacc\`a.
\newblock Explicit {B}rill-{N}oether-{P}etri general curves.
\newblock {\em Comment. Math. Helv.}, 91(3):477--491, 2016.

\bibitem[ACT21]{ACT21}
D.~Anderson, L.~Chen, and N.~Tarasca.
\newblock {K}-classes of {B}rill-{N}oether loci and a determinantal formula.
\newblock {\em International Mathematics Research Notices}, 2021.

\bibitem[AF11]{AproduFarkas11}
M.~Aprodu and G.~Farkas.
\newblock Koszul cohomology and applications to moduli.
\newblock In D.~Ellwood and E.~Previato, editors, {\em Grassmannians, moduli
  spaces, and vector bundles}. CMI/AMS, 2011.

\bibitem[ALY19]{AtanasovLarsonYang19}
A.~Atanasov, E.~Larson, and D.~Yang.
\newblock Interpolation for normal bundles of general curves.
\newblock {\em Mem. Amer. Math. Soc.}, 257(1234):v+105, 2019.

\bibitem[Ber87]{Bertram87}
A.~Bertram.
\newblock An existence theorem for {P}rym special divisors.
\newblock {\em Invent. Math.}, 90(3):669--671, 1987.

\bibitem[BF86]{BF86}
R.~Biggers and M.~Fried.
\newblock Irreducibility of moduli spaces of cyclic unramified covers of genus
  {$g$} curves.
\newblock {\em Trans. Amer. Math. Soc.}, 295(1):59--70, 1986.

\bibitem[BK96]{BallicoKeem96}
E.~Ballico and C.~Keem.
\newblock On linear series on general {$k$}-gonal projective curves.
\newblock {\em Proc. Amer. Math. Soc.}, 124(1):7--9, 1996.

\bibitem[BN74]{BN1874}
V.~Brill and M.~Noether.
\newblock Ueber die algebraischen {F}unctionen und ihre {A}nwendung in der
  {G}eometrie.
\newblock {\em Mathematische Annalen}, 7:269--310, 1874.

\bibitem[CDPR12]{tropicalBN}
F.~Cools, J.~Draisma, S.~Payne, and E.~Robeva.
\newblock A tropical proof of the {B}rill-{N}oether theorem.
\newblock {\em Adv. Math.}, 230(2):759--776, 2012.

\bibitem[CFM13]{ChenFarkasMorrison13}
D.~Chen, G.~Farkas, and I.~Morrison.
\newblock Effective divisors on moduli spaces of curves and abelian varieties.
\newblock In {\em A celebration of algebraic geometry}, volume~18 of {\em Clay
  Math. Proc.}, pages 131--169. Amer. Math. Soc., Providence, RI, 2013.

\bibitem[CJ19]{CPJ19}
K.~{Cook-Powell} and D.~Jensen.
\newblock Components of {B}rill-{N}oether loci for curves with fixed gonality.
\newblock Preprint, arXiv:1907.08366, 2019.

\bibitem[CJ20]{CPJ20}
K.~{Cook-Powell} and D.~Jensen.
\newblock Tropical methods in {H}urwitz-{B}rill-{N}oether theory.
\newblock Preprint, arXiv:2007.13877v1, 2020.

\bibitem[CJP15]{LiftingDivisors}
D.~Cartwright, D.~Jensen, and S.~Payne.
\newblock Lifting divisors on a generic chain of loops.
\newblock {\em Canad. Math. Bull.}, 58(2):250--262, 2015.

\bibitem[Cli78]{Clifford}
W.~Clifford.
\newblock On the classification of loci.
\newblock {\em Philosophical Transactions of the Royal Society of London},
  169:663--681, 1878.

\bibitem[CLRW20]{CLRW20}
S.~Creech, Y.~Len, C.~Ritter, and D.~Wu.
\newblock {P}rym-{B}rill-{N}oether loci of special curves.
\newblock Preprint, arXiv:1912.02863v2, 2020.

\bibitem[CM99]{CoppensMartens}
M.~Coppens and G.~Martens.
\newblock Linear series on a general $k$-gonal curve.
\newblock {\em Abhandlungen aus dem Mathematischen Seminar der Universit{\"a}t
  Hamburg}, 69:347, 1999.

\bibitem[CP12]{CartwrightPayne12}
D.~Cartwright and S.~Payne.
\newblock Connectivity of tropicalizations.
\newblock {\em Math. Res. Lett.}, 19(5):1089--1095, 2012.

\bibitem[CP21]{CP21}
M.~Chan and N.~Pflueger.
\newblock Euler characteristics of {B}rill-{N}oether varieties.
\newblock {\em Trans. Amer. Math. Soc.}, 374(3):1513--1533, 2021.

\bibitem[DCP95]{dCP95}
C.~De~Concini and P.~Pragacz.
\newblock On the class of {B}rill-{N}oether loci for {P}rym varieties.
\newblock {\em Math. Ann.}, 302(4):687--697, 1995.

\bibitem[Deb00]{Debarre00}
O.~Debarre.
\newblock Th\'{e}or\`emes de {L}efschetz pour les lieux de
  d\'{e}g\'{e}n\'{e}rescence.
\newblock {\em Bull. Soc. Math. France}, 128(2):283--308, 2000.

\bibitem[EH83a]{EisenbudHarris83b}
D.~Eisenbud and J.~Harris.
\newblock Divisors on general curves and cuspidal rational curves.
\newblock {\em Invent. Math.}, 74(3):371--418, 1983.

\bibitem[EH83b]{EisenbudHarris83c}
D.~Eisenbud and J.~Harris.
\newblock A simpler proof of the {G}ieseker-{P}etri theorem on special
  divisors.
\newblock {\em Invent. Math.}, 74(2):269--280, 1983.

\bibitem[EH87a]{EisenbudHarris87}
D.~Eisenbud and J.~Harris.
\newblock Irreducibility and monodromy of some families of linear series.
\newblock {\em Ann. Sci. \'Ecole Norm. Sup. (4)}, 20(1):65--87, 1987.

\bibitem[EH87b]{EisenbudHarris87b}
D.~Eisenbud and J.~Harris.
\newblock The {K}odaira dimension of the moduli space of curves of genus {$\geq
  23$}.
\newblock {\em Invent. Math.}, 90(2):359--387, 1987.

\bibitem[EH89]{EisenbudHarris89b}
D.~Eisenbud and J.~Harris.
\newblock Progress in the theory of complex algebraic curves.
\newblock {\em Bull. Amer. Math. Soc. (N.S.)}, 21(2):205--232, 1989.

\bibitem[Far09]{Farkas09b}
G.~Farkas.
\newblock Koszul divisors on moduli spaces of curves.
\newblock {\em Amer. J. Math.}, 131(3):819--867, 2009.

\bibitem[FJP20]{M23}
G.~Farkas, D.~Jensen, and S.~Payne.
\newblock The {K}odaira dimensions of ${M}_{22}$ and ${M}_{23}$.
\newblock Preprint, arXiv:2005.00622v1, 2020.

\bibitem[FJP21]{R13}
G.~Farkas, D.~Jensen, and S.~Payne.
\newblock The non-abelian {B}rill-{N}oether divisor on
  {$\overline{\mathcal{M}}_{13}$} and the {K}odaira dimension of
  {$\overline{\mathcal{R}}_{13}$}.
\newblock Preprint, arXiv:2110.09553, 2021.

\bibitem[FL81]{FultonLazarsfeld81}
W.~Fulton and R.~Lazarsfeld.
\newblock On the connectedness of degeneracy loci and special divisors.
\newblock {\em Acta Math.}, 146(3-4):271--283, 1981.

\bibitem[FP05]{FarkasPopa05}
G.~Farkas and M.~Popa.
\newblock Effective divisors on {$\overline{\mathcal M}_g$}, curves on {$K3$}
  surfaces, and the slope conjecture.
\newblock {\em J. Algebraic Geom.}, 14(2):241--267, 2005.

\bibitem[Ful69]{Fulton69}
W.~Fulton.
\newblock Hurwitz schemes and irreducibility of moduli of algebraic curves.
\newblock {\em Ann. of Math. (2)}, 90:542--575, 1969.

\bibitem[GH80]{GriffithsHarris80}
P.~Griffiths and J.~Harris.
\newblock On the variety of special linear systems on a general algebraic
  curve.
\newblock {\em Duke Math. J.}, 47(1):233--272, 1980.

\bibitem[Gie82]{Gieseker82}
D.~Gieseker.
\newblock Stable curves and special divisors: {P}etri's conjecture.
\newblock {\em Invent. Math.}, 66(2):251--275, 1982.

\bibitem[Har82]{Harris82}
J.~Harris.
\newblock {\em Curves in projective space}, volume~85 of {\em S\'eminaire de
  Math\'ematiques Sup\'erieures}.
\newblock Presses de l'Universit\'e de Montr\'eal, Montreal, Que., 1982.
\newblock With the collaboration of David Eisenbud.

\bibitem[Har84]{Harris84}
J.~Harris.
\newblock On the {K}odaira dimension of the moduli space of curves. {II}. {T}he
  even-genus case.
\newblock {\em Invent. Math.}, 75(3):437--466, 1984.

\bibitem[Hir85]{Hirschowitz85}
A.~Hirschowitz.
\newblock La m\'{e}thode d'{H}orace pour l'interpolation \`a plusieurs
  variables.
\newblock {\em Manuscripta Math.}, 50:337--388, 1985.

\bibitem[HM82]{HarrisMumford82}
J.~Harris and D.~Mumford.
\newblock On the {K}odaira dimension of the moduli space of curves.
\newblock {\em Invent. Math.}, 67(1):23--88, 1982.
\newblock With an appendix by William Fulton.

\bibitem[JP14]{tropicalGP}
D.~Jensen and S.~Payne.
\newblock Tropical independence {I}: {S}hapes of divisors and a proof of the
  {G}ieseker-{P}etri theorem.
\newblock {\em Algebra Number Theory}, 8(9):2043--2066, 2014.

\bibitem[JR21]{JensenRanganathan}
D.~Jensen and D.~Ranganathan.
\newblock Brill-{N}oether theory for curves of a fixed gonality.
\newblock {\em Forum Math. Pi}, 9:e1, 33, 2021.

\bibitem[Kem71]{Kempf71}
G.~Kempf.
\newblock {\em Schubert methods with an application to algebraic curves}.
\newblock Publ. Math. Centrum, 1971.

\bibitem[KL72]{KleimanLaksov72}
S.~Kleiman and D.~Laksov.
\newblock On the existence of special divisors.
\newblock {\em Amer. J. Math.}, 94:431--436, 1972.

\bibitem[Lar17]{Larson17}
E.~Larson.
\newblock The maximal rank conjecture.
\newblock Preprint, arXiv:1711.04906, 2017.

\bibitem[Lar18]{Larson18}
E.~Larson.
\newblock Degenerations of curves in projective space and the maximal rank
  conjecture.
\newblock Preprint, arXiv:1809.05980, 2018.

\bibitem[Lar20]{Larson20}
E.~Larson.
\newblock The maximal rank conjecture for sections of curves.
\newblock {\em J. Algebra}, 555:223--245, 2020.

\bibitem[Lar21]{Larson21}
H.~Larson.
\newblock A refined {B}rill-{N}oether theory over {H}urwitz spaces.
\newblock {\em Invent. Math.}, 224(3):767--790, 2021.

\bibitem[Laz86]{Lazarsfeld86}
R.~Lazarsfeld.
\newblock Brill-{N}oether-{P}etri without degenerations.
\newblock {\em J. Differential Geom.}, 23(3):299--307, 1986.

\bibitem[LLM{\etalchar{+}}14]{LLMSSZ}
T.~Lam, L.~Lapointe, J.~Morse, A.~Schilling, M.~Shimozono, and M.~Zabrocki.
\newblock {\em {$k$}-{S}chur functions and affine {S}chubert calculus},
  volume~33 of {\em Fields Institute Monographs}.
\newblock Springer, New York; Fields Institute for Research in Mathematical
  Sciences, Toronto, ON, 2014.

\bibitem[LLV20]{LLV}
E.~Larson, H.~Larson, and I.~Vogt.
\newblock Global {B}rill-{N}oether theory over the {H}urwitz space.
\newblock Preprint, arXiv:2008.10765, 2020.

\bibitem[LU20]{LenUlirsch21}
Y.~Len and M.~Ulirsch.
\newblock Skeletons of {P}rym varieties and {B}rill-{N}oether theory.
\newblock Preprint, arXiv:1902.09410v2, 2020.

\bibitem[Mar46]{Maroni}
A.~Maroni.
\newblock Le serie lineari speciali sulle curve trigonali.
\newblock {\em Ann. Mat. Pura Appl. (4)}, 25:343--354, 1946.

\bibitem[Mar96]{Martens96}
G.~Martens.
\newblock On curves of odd gonality.
\newblock {\em Arch. Math. (Basel)}, 67(1):80--88, 1996.

\bibitem[Oss14]{Osserman14}
B.~Osserman.
\newblock A simple characteristic-free proof of the {B}rill-{N}oether theorem.
\newblock {\em Bull. Braz. Math. Soc. (N.S.)}, 45(4):807--818, 2014.

\bibitem[Pfl17]{Pflueger17b}
N.~Pflueger.
\newblock Brill-{N}oether varieties of $k$-gonal curves.
\newblock {\em Adv. Math.}, 312:46--63, 2017.

\bibitem[Ste98]{Steffen98}
F.~Steffen.
\newblock A generalized principal ideal theorem with an application to
  {B}rill-{N}oether theory.
\newblock {\em Invent. Math.}, 132(1):73--89, 1998.

\bibitem[TiB03]{TiB03}
M.~Teixidor~i {B}igas.
\newblock Injectivity of the symmetric map for line bundles.
\newblock {\em Manuscripta Math.}, 112(4):511--517, 2003.

\bibitem[Wel85]{Welters85}
G.~Welters.
\newblock A theorem of {G}ieseker-{P}etri type for {P}rym varieties.
\newblock {\em Ann. Sci. \'{E}cole Norm. Sup. (4)}, 18(4):671--683, 1985.

\end{thebibliography}

\end{document}